\documentclass[12pt,reqno]{amsart}
\usepackage{amssymb,enumerate}

\usepackage{pdfpages}

\DeclareMathAlphabet{\E}{U}{eus}{m}{n}     % Euler script

\newcommand{\g}{{\gamma}}      

\newcommand{\Z}{{\E{Z}}}

\newcommand{\V}{{\mathcal V}}

\newcommand{\A}{{\mathcal A}}

\newcommand{\PP}{{\mathbb P}}

\newcommand{\kk}{{\Bbbk}}

\renewcommand{\L}{{\mathfrak L}}

\newcommand{\p}{\mathfrak p}

\newtheorem{thm}{Theorem}[section]
\newtheorem{lemma}[thm]{Lemma}
\newtheorem{cor}[thm]{Corollary}
\newtheorem{prop}[thm]{Proposition}

\theoremstyle{definition}
\newtheorem{conj}[thm]{Conjecture}
\newtheorem{defn}[thm]{Definition}

\newtheorem{rmk}[thm]{Remark}

\newtheorem{ack}{Acknowledgments\!\!}  

\newtheorem{Pf}{Proof$\!\!$}         
         
\newenvironment{pf}{\begin{Pf}}{\qed\end{Pf}}

\DeclareMathSymbol{\twoheadrightarrow}  {\mathrel}{AMSa}{"10}
\newcounter{letter}
\renewcommand{\theletter}{\rom{(}\alph{letter}\rom{)}}

\newcounter{rnum}
\renewcommand{\thernum}{\rom{(}\roman{rnum}\rom{)}}

\begin{document}
\baselineskip21pt

%\hspace*{\fill}\today

\title[Line Scheme of a Certain Family of Quantum $\PP^3$s]%
{The One-Dimensional Line Scheme of a\\[2mm]
Certain Family of Quantum $\PP^3$s}

%\subjclass{14A22, 16S37, 16S38}%
%14A22 Noncommutative algebraic geometry
%16S37 Quadratic and Koszul algebras
%16S38 Rings arising from non-commutative algebraic geometry

\keywords{line scheme, point scheme, elliptic curve, % 
regular algebra, Pl\"ucker coordinates.%
%\newline
%\indent This work was supported in part by NSF grants DMS-0900239 and
%DMS-1302050
}
%\rule[-5mm]{0cm}{0cm}}%

\maketitle

\vspace*{0.1in}

\baselineskip15pt

\renewcommand{\thefootnote}{\fnsymbol{footnote}}
%\centerline{\sc Michaela Vancliff\footnote{The first author was
%supported in part by NSF grants DMS-0457022 \& DMS-0900239.\\[-3mm]}}
%\centerline{Department of Mathematics, P.O.~Box 19408}
%\centerline{University of Texas at Arlington,
%Arlington, TX 76019-0408}
%\centerline{{\sf vancliff@uta.edu}}
%\centerline{{\sf www.uta.edu/math/vancliff}}

\begin{center}
\begin{tabular}[t]{c}
\textsc{Richard G.\ Chandler$^*$}\\
richard.chandler@mavs.uta.edu\\
students.uta.edu/rg/rgc7061
\end{tabular}
\quad  and \quad  
\begin{tabular}[t]{c}
\textsc{Michaela Vancliff}\\
vancliff@uta.edu\\
www.uta.edu/math/vancliff
\end{tabular}%
\hspace*{-6.5mm}\footnote{%
This work was supported in part by NSF grants DMS-0900239 and
DMS-1302050.}\\[3mm]
Department of Mathematics, P.O.~Box 19408\\
University of Texas at Arlington,\\
Arlington, TX 76019-0408
\end{center}

\setcounter{page}{1}
\thispagestyle{empty}

\bigskip
\bigskip

\begin{abstract}
\baselineskip15pt
A quantum~$\PP^3$ is a noncommutative analogue of a polynomial ring on four 
variables, and, herein, it is taken  to be a regular algebra of global 
dimension four.  It is well known that if a generic quadratic 
quantum~$\PP^3$ exists, then it has a point scheme consisting of exactly 
twenty
distinct points and a one-dimensional line scheme. In this article, we
compute the line scheme of a family of algebras whose generic member is a 
candidate for a generic quadratic quantum~$\PP^3$. We find that, as
a closed subscheme of $\PP^5$, the line 
scheme of the generic member is the union of seven curves; namely, 
a nonplanar elliptic curve in a $\PP^3$, four planar elliptic curves and 
two nonsingular conics. 
\end{abstract}

\baselineskip18pt

%\newpage
\bigskip
\bigskip

%%%%%%%%%%%%%%%%%%%%%%%%%%%%%%%%%%%%%%%%%%%%%%%%%%%%%%%%%%%%%%%%%%%%%%
% INTRO
%%%%%%%%%%%%%%%%%%%%%%%%%%%%%%%%%%%%%%%%%%%%%%%%%%%%%%%%%%%%%%%%%%%%%%

\section*{Introduction}

A regular algebra of global dimension~$n$ is often viewed as a
noncommutative analogue of a polynomial ring on $n$~variables.
Generalizing the language in \cite{A}, such an algebra is sometimes called 
a quantum~$\PP^{n-1}$. In \cite{ATV1}, quantum~$\PP^2$s were classified
according to their point schemes, with the point scheme of the most generic
quadratic quantum~$\PP^2$ depicted by an elliptic curve in $\PP^2$.

Consequently, a similar description is desired for quadratic 
quantum~$\PP^3$s using their point schemes or their line schemes, where the
definition of line scheme was given in \cite{SV1}.  However, to date, very 
few line schemes of quadratic quantum~$\PP^3$s are known, especially of 
algebras that are candidates for
generic quadratic quantum~$\PP^3$s. As explained in \cite{Vmsri}, if a 
generic quadratic quantum~$\PP^3$ exists, then it has a point scheme 
consisting of exactly twenty distinct points and a one-dimensional line 
scheme.  Hence, in this article, we compute the line scheme of a family of 
algebras that appeared in \cite[\S5]{CV}, and whose generic member is a 
candidate for a generic quadratic quantum~$\PP^3$.  

The article is outlined as follows. Section~1 begins with some definitions,
including the introduction of the family of algebras considered herein.
The point schemes of the algebras are computed in Section~2 in
Proposition~\ref{ptscheme}, whereas 
Sections~3 and 4 are devoted to the computation of the line scheme and 
identifying the lines in $\PP^3$ to which the points of the line scheme 
correspond. In particular, our main results are Theorems~\ref{L'g},  
\ref{Lg} and \ref{6lines}.  In the first two, we prove that the line scheme 
of the generic member is the union of seven curves; namely, a nonplanar 
elliptic curve in a $\PP^3$ (a spatial elliptic curve), four planar 
elliptic curves and two nonsingular conics. 
In Theorem~\ref{6lines}, we find that if $p$ is one of the generic points 
of the point scheme, then there are exactly six distinct lines of the line 
scheme that pass through $p$. 
An Appendix is provided in Section~\ref{app} that lists polynomials that 
are used throughout the article. 

It is hoped that data 
from the one-dimensional line scheme of any potentially generic quadratic 
quantum~$\PP^3$ will motivate conjectures and future research in the
subject. In fact, the results herein suggest that the line scheme of 
the most generic quadratic quantum~$\PP^3$ is conceivably the union of 
two spatial elliptic curves and four planar elliptic curves (see 
Conjecture~\ref{conj}).

%%%%%%%%%%%%%%%%%%%%%%%%%%%%%%%%%%%%%%%%%%%%%%%%%%%%%%%%%%%%%%%%%%%%%
% END  INTRODUCTION
% START  SECTION  1
%%%%%%%%%%%%%%%%%%%%%%%%%%%%%%%%%%%%%%%%%%%%%%%%%%%%%%%%%%%%%%%%%%%%%

\bigskip
%\bigskip

\section{The Algebras} \label{sec1}

In this section, we introduce the algebras from \cite[\S5]{CV} that are
considered in this article.

Throughout the article, $\kk$ denotes an algebraically closed field and 
$M(n,\, \kk)$ denotes the vector space of $n \times n$ matrices with entries
in~$\kk$.  If $V$ is a vector space, then $V^\times$ will denote the nonzero
elements in $V$, and $V^{\text{*}}$ will denote the vector-space dual of~$V$.
In this section, we take char$(\kk) \neq 2$, but, in Sections~\ref{sec3}
and \ref{sec4}, we assume char$(\kk) = 0$ owing to the computations in 
those sections.

\begin{defn}\label{Ag}\cite[\S5]{CV}
Let $\g \in \kk^\times$ and write $\A(\g)$ for the $\kk$-algebra
on generators $x_1, \ldots, x_4$ with defining relations:
\[
\begin{array}{lll}
x_4 x_1 = i x_1 x_4, \qquad & x_3^2 = x_1^2, \qquad &
x_3 x_1 = x_1 x_3 - x_2^2,\\[3mm]
x_3 x_2 = i x_2 x_3, & x_4^2 = x_2^2, & x_4 x_2 = x_2 x_4 - \g
x_1^2,
\end{array}
\]
where $i^2 = -1$.
\end{defn}

By construction of $\A(\g)$ in \cite{CV}, $\A(\g)$ is a regular noetherian 
domain of global dimension four with Hilbert series the same as that of the
polynomial ring on four variables.  As remarked 
in \cite{CV}, the special 
member $\A(1)$ was studied in \cite{ST} and, if $\g^2 \neq 4$, then 
$\A(\g)$ has a finite point scheme consisting of twenty distinct points and 
a one-dimensional line scheme.  Since the computation of the point scheme 
was omitted from \cite{CV}, we will outline the computation of it in 
Section~\ref{sec2}.

It should be noted that $\A(\g) \cong \A(-\g)$, for all $\g\in \kk^\times$,
under the map 
that sends $x_2 \mapsto -x_2$ and $x_k \mapsto x_k$ for all $k \neq 2$. 
There also exist antiautomorphisms of $\A(\g)$ defined by 
\[
\begin{array}{l}
\psi_1: \ 
x_1 \leftrightarrow x_3 \text{ \ and \ } x_2 \leftrightarrow x_4,\\[2mm]
\psi_2: \ 
x_2 \leftrightarrow \lambda x_3 \text{ \ and \ } x_4 \leftrightarrow
\lambda x_1,
\end{array}
\]
where $\lambda \in \kk^\times$ with $\lambda^4 = \g$. These latter maps
will be useful in Sections~\ref{sec3} and \ref{sec4}.

The reader should note that the point scheme given in \cite{ST} for
$\A(1)$ has some sign errors in the formulae. Moreover, $\A(1)$ was studied 
in \cite{G} in the context of finding the scheme of lines associated
to each point of the point scheme.

For background material on point modules, line modules, point schemes,
line schemes, regular algebras and some of the historical development of 
the subject, the reader is referred to \cite{Vmsri}.

%%%%%%%%%%%%%%%%%%%%%%%%%%%%%%%%%%%%%%%%%%%%%%%%%%%%%%%%%%%%%%%%%%%%%
% END  SECTION  1
% START  SECTION  2
%%%%%%%%%%%%%%%%%%%%%%%%%%%%%%%%%%%%%%%%%%%%%%%%%%%%%%%%%%%%%%%%%%%%%

\bigskip
%\newpage

\section{The Point Scheme of $\A(\g)$}\label{sec2}

In this section, we compute the point scheme of the algebras $\A(\g)$
given in Definition~\ref{Ag}. Our method follows that of \cite{ATV1},
and we continue to assume that char$(\kk) \neq 2$ in this section. 

Let $V = \sum_{i=1}^4 \kk x_i$. Following \cite{ATV1},
we write the relations of $\A(\g)$ in the form $Mx = 0$, where $M$ is a 
$6\times 4$ matrix and $x$ is the column vector given by $x^T = (x_1,
\ldots, x_4)$. Thus, we may take $M$ to be the matrix 
\[
M= \begin{bmatrix} 
x_4 & 0 & 0 & -i x_1\\
0 & x_3 & -i x_2 & 0\\
x_1 & 0 & -x_3 & 0\\
0 & x_2 & 0 & -x_4\\
x_3 & x_2 & -x_1 & 0\\
\g x_1 & x_4 & 0 & -x_2
\end{bmatrix}
,
\]
and, by \cite{ATV1}, the point scheme of $\A(\g)$ can be identified with 
the zero locus, $\p(\g)$, in $\PP(V^{\text{*}})$ of all the 
$4 \times 4$~minors of $M$.  Fifteen polynomials given by these minors are 
listed in Section~\ref{app1} in the Appendix.
We will prove that, if $\g^2 \neq 4$, then $\p(\g)$ is finite with twenty 
distinct points. 

Let $p = (\alpha_1, \ldots , \alpha_4) \in \p(\g)$. If $\alpha_1 = 0$,
then it is straightforward to prove that $p$ is one of the points 
$e_2 = (0,\,1,\,0,\,0)$, $e_3 = (0,\,0,\,1,\,0)$, $e_4 = (0,\,0,\,0,\,1)$. 
Thus, we assume $\alpha_1 = 1$. If, in addition, $\alpha_4 = 0$, then
rank$(M) = 0$ if and only if $\alpha_2 = 0 = \alpha_3$, so we obtain
the point $e_1 = (1,\,0,\,0,\,0)$. Hence, we may assume $\alpha_1 = 1$
and $\alpha_4 \neq 0$. 

With this assumption, a computer-algebra program such as Wolfram's
Mathematica yields three polynomials that determine the remaining
closed points in $\p(\g)$:
\[
\rho_1 = x_4^8 -4 x_4^4 + \g^2,\quad
\rho_2 = x_3^2 - i x_3 x_4^2 - 1,\quad
\rho_3 = \g x_2 - 2i x_4^3 + x_3 x_4^5.
\]
(In fact, \ref{app1}.\ref{pp1}, \ref{app1}.\ref{pp2} and 
\ref{app1}.\ref{pp5} evaluated at $x_1 = 1$ generate the other polynomials 
in Section~\ref{app1} evaluated at $x_1 = 1$, and determine 
$\rho_1,\, \rho_2,\, \rho_3$.)
Since $\rho_1 = 0$ if and only if $(x_4^4 -2)^2 = 4 - \g^2$, we
find that $\rho_1$ has eight distinct zeros if and only if
$\g^2 \neq 4$; if $\g^2 = 4$, then $\rho_1$ has exactly four distinct
zeros, each of multiplicity two. 
Given a zero $x_4$ to $\rho_1$, the equation~$\rho_2 = 0$ has 
a unique solution for $x_3$ if and only if $x_4^4 = 4$, but this implies 
$\rho_1 \neq 0$ as $\g \neq 0$, which is false; hence $\rho_2$ has two 
distinct zeros for all $\g \in \kk^\times$.

The following remark will be useful in the proof of 
Proposition~\ref{ptscheme}.
\begin{rmk} \label{20pts} (cf., \cite{Vmsri})
If the zero locus $\mathfrak z$ of the defining relations of a quadratic 
algebra on four generators with six defining relations is finite, then 
$\mathfrak z$ consists of twenty points counted with multiplicity.
\end{rmk}

\begin{prop}\label{ptscheme}
Let $\A(\g)$ and $\p(\g)$ be as above and let $\Z_\g$ denote the scheme of 
zeros of $\rho_1,\, \rho_2,\, \rho_3$ in $\PP(V^{\text{\em *}})$.
\begin{enumerate}
\itemsep1mm
\item[\textup{(a)}]
           For every $\g \in \kk^\times$,
	   $\p(\g) = \{ e_1, \ldots, e_4 \} \cup \Z_\g$.
\item[\textup{(b)}]
          If $\g^2 \neq 4$, then $\p(\g)$ has exactly twenty distinct points.
\item[\textup{(c)}] 
           If $\g^2 = 4$, then $\p(\g)$ has exactly twelve distinct 
           points; the eight closed points of $\Z_\g$ have multiplicity two 
	   in $\p(\g)$ and the remaining four points of $\p(\g)$
	   each have multiplicity one.
\item[\textup{(d)}] 
           For every $\g \in \kk^\times$, the closed points in 
           $\PP(V^{\text{\em *}}) \times \PP(V^{\text{\em *}})$ on which 
	   the defining relations of~$\A(\g)$ vanish are given by: 
	   $(e_1,\, e_2)$, $(e_2,\, e_1)$, $(e_3,\, e_4)$, $(e_4,\, e_3)$ 
           and points of the form
           \[\big((1,\, \alpha_2,\, \alpha_3,\, \alpha_4), \   
           (1,\, i\alpha_2\alpha_3^{-2},\, \alpha_3^{-1},\, -i\alpha_4)
	   \big),\] 
           where $(1,\, \alpha_2,\, \alpha_3,\, \alpha_4) \in \Z_\g$ and 
	   $i^2 = -1$.
\end{enumerate}
\end{prop}
\begin{pf}
The preceding discussion proves that if $\g^2 \neq 4$, then the number of 
distinct closed points in $\p(\g)$ is twenty, so, by
Remark~\ref{20pts}, (b) follows.  On the other hand, if $\g^2 = 4$, then 
the zeros of $\rho_1$ have multiplicity two, so, counting multiplicity, 
the eight distinct points in $\Z_\g$ have multiplicity two. Thus, each $e_i$
has multiplicity one, by Remark~\ref{20pts}.  Hence, (c) and (a) follow.
Part (d) is easily verified by computation with the matrix $M$ using
polynomials 
\ref{app1}.\ref{pp1}, \ref{app1}.\ref{pp2} and \ref{app1}.\ref{pp10} in
the Appendix.
\end{pf}

\begin{cor}
For all $\g \in \kk^\times$, there exists an automorphism 
$\sigma : \p(\g) \to \p(\g)$ which, on closed points, is defined by: 
\[
\begin{array}{c}
e_1 \leftrightarrow e_2,\quad e_3 \leftrightarrow e_4,\\[2mm] 
\sigma\big((1,\, \alpha_2,\, \alpha_3,\, \alpha_4) \big) = 
(1,\, i\alpha_2\alpha_3^{-2},\, \alpha_3^{-1},\, -i\alpha_4)\\[2mm]
\end{array}
\]
for all $(1, \,\alpha_2,\, \alpha_3,\, \alpha_4) \in \Z_\g$. Hence,
on the closed points of $\p(\g)$, $\sigma$ has two orbits of length two
and $n$ orbits of length four, where $n = 4$ if \;$|\Z_\g| = 16$ and 
$n = 2$ if \;$|\Z_\g| = 8$.
\end{cor}
\begin{pf}
The fact the map exists on the closed points of $\p(\g)$ is a
consequence of Proposition~\ref{ptscheme}(d); its existence on the scheme
follows from \cite[Theorem~4.1.3]{LSV}. The size of the orbits may be
verified by computation. 
\end{pf}

%%%%%%%%%%%%%%%%%%%%%%%%%%%%%%%%%%%%%%%%%%%%%%%%%%%%%%%%%%%%%%%%%%%%%
% END  SECTION  1
% START  SECTION  2
%%%%%%%%%%%%%%%%%%%%%%%%%%%%%%%%%%%%%%%%%%%%%%%%%%%%%%%%%%%%%%%%%%%%%

\bigskip
%\newpage

\section{The Line Scheme of $\A(\g)$}\label{sec3}

In this section, we compute the line scheme $\L(\g)$ of the algebras 
$\A(\g)$ as a closed subscheme of~$\PP^5$.
Our arguments follow the method given in \cite{SV2}, which is summarized  
below
in Section~3.1. In Section~3.2, we compute the closed points of the line 
scheme, and, in Section~3.3, we prove that the line scheme is a reduced 
scheme, and so is given by its closed points. The main results of this
section are Theorems~\ref{L'g} and \ref{Lg}.
Henceforth, we assume that char$(\kk) = 0$.

\medskip

\subsection{Method}\label{method}\hfill

In \cite{SV2}, a method was given for computing the line scheme of any 
quadratic algebra on four generators that is a domain and has Hilbert 
series the same as that of the polynomial ring on four variables. In
this subsection, we summarize that method while applying it to $\A(\g)$; 
further details may be found in \cite{SV2}.

The first step in the process is to compute the Koszul dual of $\A(\g)$. 
This produces a quadratic algebra on four generators with ten defining 
relations. One then rewrites those ten relations in the form of a matrix 
equation similar to that used in Section~\ref{sec2}; in this case,
however, it yields the equation $\hat M z = 0$, where 
$z^T = (z_1, \ldots , z_4)$ (where $\{ z_1, \ldots , z_4\}$ is the dual 
basis in $V^{\text{*}}$ to $\{x_1, \ldots , x_4\}$) and $\hat M$ is a
$10\times 4$ matrix whose entries are linear forms in the $z_i$. 

One then produces a $10 \times 8$ matrix from $\hat M$ by concatenating two 
$10\times 4$ matrices, the first of which is obtained from $\hat M$ by 
replacing every $z_i$ in $\hat M$ by $u_i\in\kk$, and the second is obtained
from $\hat M$ by replacing every $z_i$ in $\hat M$ by $v_i\in\kk$,
where $(u_1, \ldots, u_4 ), \  (v_1, \ldots, v_4 ) \in \PP^3$.  
For $\A(\g)$, this process yields the following $10 \times 8$ matrix:
\[
\small
\mathcal M(\g) = 
\begin{bmatrix}
0 & u_1 & 0 & 0 & 0 & v_1 & 0 & 0\\
u_2 & 0 & 0 & 0 & v_2 & 0 & 0 & 0\\
0 & 0 & 0 & u_3 & 0 & 0 & 0 & v_3\\
0 & 0 & u_4 & 0 & 0 & 0 & v_4 & 0\\
u_3 & 0 & u_1 & 0 & v_3 & 0 & v_1 & 0\\
0 & u_4 & 0 & u_2 & 0 & v_4 & 0 & v_2\\
-u_4 & 0 & 0 & i u_1 & -v_4 & 0 & 0 & i v_1\\
0 & -u_3 & i u_2 & 0 & 0 & -v_3 & i v_2 & 0\\
u_1 & 0 & u_3 & \g u_2 & v_1 & 0 & v_3 & \g v_2\\
0 & u_2 & u_1 & u_4 & 0 & v_2 & v_1 & v_4
\end{bmatrix}.
\]

\quad

Each of the forty-five $8\times 8$ minors of $\mathcal M(\g)$ is a 
bihomogeneous 
polynomial of bidegree $(4, 4)$ in the $u_i$ and $v_i$, and so each such 
minor is a linear combination of products of polynomials of the form 
$N_{ij}=u_i v_j- u_j v_i$, where $1 \leq i < j\leq 4$. 
Hence, $\mathcal M(\g)$ yields forty-five quartic polynomials in the six 
variables~$N_{ij}$.  Following \cite{SV2}, one then applies the map:\\[-2mm]
\[
\begin{array}{l}
N_{12}\mapsto M_{34}, \quad N_{13}\mapsto -M_{24}, \quad N_{14}\mapsto M_{23},\\[3mm]
N_{23}\mapsto M_{14}, \quad N_{24}\mapsto -M_{13}, \quad N_{34}\mapsto M_{12},
\end{array}
\]
\quad\\[-1mm]
to the polynomials, which yields forty-five quartic polynomials in the 
Pl\"ucker coordinates~$M_{ij}$ on~$\PP^5$. 

The line scheme $\L(\g)$ of $\A(\g)$ may be realised in~$\PP^5$ as the 
scheme of zeros of these forty-five polynomials in the $M_{ij}$ together 
with the Pl\"ucker 
polynomial $P = M_{12}M_{34}-M_{13}M_{24}+M_{14}M_{23}$. For $\A(\g)$,
these polynomials were found by using Wolfram's Mathematica and are listed 
in Section~\ref{app2} of the Appendix.  

In the remainder of this section, we compute and describe $\L(\g)$
as a subscheme of $\PP^5$.  The lines in $\PP(V^{\text{*}})$ that 
correspond to the points of $\L(\g)$ are described in Section~\ref{sec4}.

%\bigskip
%\newpage

\subsection{Computing the Closed Points of the Line Scheme}\hfill

Our procedure in this subsection focuses on finding the closed points of 
the line scheme $\L(\g)$ of $\A(\g)$; in the next subsection, we will prove 
that $\L(\g)$ is reduced and so is given by its closed points. We
denote the variety of closed points of $\L(\g)$  by $\L'(\g)$ and the
zero locus of a set $S$ of polynomials by $\V(S)$.

Subtracting the polynomials \ref{app2}.\ref{p1} and \ref{app2}.\ref{p2} 
produces 
$M_{14}M_{23}M_{24}^2$.  If $M_{14} = M_{23} = M_{24} = 0$, then
$M_{12} = 0 = M_{34}$, so there is a unique solution in this case. This
leaves six cases to consider:\\[-2mm]
\begin{center}
\begin{tabular}{rlrl}
(I)& $M_{14}M_{23} \neq 0,\ M_{24}=0$,\hspace*{9mm} &
(IV) & $M_{23} \neq 0,\ M_{14}=0=M_{24}$,\\[2mm]
(II)& $M_{23}M_{24} \neq 0,\ M_{14}=0$, &
(V)& $M_{14} \neq 0,\ M_{23}=0=M_{24}$, \\[2mm]
(III) & $M_{14}M_{24} \neq 0,\ M_{23}=0$, &
(VI) & $M_{24} \neq 0,\ M_{14}=0=M_{23}$.
\end{tabular}
\end{center}

\bigskip

We will outline the analysis for (I), (II), (IV) and (VI); the other cases 
follow from these four cases by using the map $\psi_1$ defined in 
Section~\ref{sec1}. 
In applying the map $\psi_1$, the reader should recall that 
$M_{ji} = - M_{ij}$ for all $i \neq j$.

\bigskip

\noindent
\textbf{Case (I):}
$M_{14}M_{23} \neq 0$ and $M_{24}=0$.\\
With the assumption that $M_{24} = 0$, a computation of a Gr\"obner 
basis yields several polynomials, one of which is 
$M_{13}^2 M_{14} M_{23}$.  Hence, $M_{13} = 0$, and another 
computation of a Gr\"obner basis yields several polynomials, two of
which are: 
\[
\begin{array}{c}
M_{14} M_{23}+M_{12} M_{34},\\[3mm]

M_{34}^4-M_{14}^2 M_{34}^2-M_{23}^2 M_{34}^2+\g M_{14} M_{23}
M_{34}^2+M_{14}^2 M_{23}^2,
\end{array}
\]
so that, in particular, $M_{12} M_{34} \neq 0$.
Using the first polynomial to substitute for $M_{14} M_{23}$, and using
the assumption that $M_{34} \neq 0$, we find that the second polynomial 
vanishes if and only if 
$M_{12}^2 + M_{34}^2 + \g M_{14} M_{23} - M_{14}^2 - M_{23}^2 = 0$.
Another computation of a Gr\"obner basis yields only these polynomials,
so that this case provides the component 
\[
\L_1 = \V(\, M_{13}, \ M_{24}, \ M_{14} M_{23}+M_{12} M_{34}, \ 
M_{12}^2 + M_{34}^2 + \g M_{14} M_{23} - M_{14}^2 - M_{23}^2
\, ).
\]
In Theorem~\ref{L'g}, we will prove that $\L_1$  is irreducible if and 
only if $\g^2 \neq 16$. Here we show that if $\g^2 = 16$, then $\L_1$ is the 
union of two nonsingular conics. Since $\A(4) \cong \A(-4)$, it suffices to 
consider $\g = 4$. In fact, let $\alpha \in \kk$ and let 
\[Q = M_{12}^2 + M_{34}^2 + \g M_{14} M_{23} - M_{14}^2 - M_{23}^2 +
2\alpha ( M_{14} M_{23}+M_{12} M_{34} ), 
\] 
and associate to $Q$ the symmetric matrix 
\[
\begin{bmatrix}
1 & 0 & 0 & \alpha \\[1mm]
0 & -1 & \alpha+\frac{\g}{2} & 0 \\[1mm]
0 & \alpha+\frac{\g}{2} & -1 & 0 \\[1mm]
\alpha & 0 & 0 & 1
\end{bmatrix},
\]
which has rank at most two if and only if $Q$ factors. This happens if
and only if $(\g,\, \alpha ) = (\pm4,\, \mp 1)$.  It follows that if  
$\g=4$, then
\[
Q = (M_{12} - M_{34} + M_{14} - M_{23}) (M_{12} - M_{34} - M_{14} + M_{23}),
\]
and $\L_1= \L_{1a} \cup \L_{1b}$, where 
\[
\begin{array}{c}
\L_{1a} = \V(\, M_{13}, \ M_{24}, \ M_{14} M_{23}+M_{12} M_{34},  \ 
             M_{12} + M_{14} - M_{23} - M_{34}\,),\\[2mm]
\L_{1b} = \V(\, M_{13}, \ M_{24}, \ M_{14} M_{23}+M_{12} M_{34},  \ 
             M_{12} - M_{14} + M_{23} - M_{34}\,),
\end{array}
\]
and each of $\L_{1a}$ and $\L_{1b}$ is a nonsingular conic,  since
using the last polynomial in each case to 
substitute for $M_{12}$ in $M_{14} M_{23}+M_{12} M_{34}$ yields a 
rank-3 quadratic form in each case.
Moreover, $\L_{1b}$ is $\psi_1$ applied to $\L_{1a}$.

\bigskip

\noindent
\textbf{Case (II):}
$M_{23}M_{24} \neq 0$ and $M_{14}=0$.\\
With the assumption that $M_{14} = 0$, a computation of a Gr\"obner 
basis yields several polynomials, two of which are 
$M_{13} M_{23} M_{24}^2$ and $M_{23} M_{24} M_{34}^2$.  
Hence, $M_{13} = M_{34} = 0$. With these additional criteria, another 
computation of a Gr\"obner basis yields exactly three polynomials: 
$M_{12}f$, $M_{23}f$, $M_{24}f$, where 
$f = M_{12}^3 - M_{12} M_{23}^2 -i M_{23} M_{24}^2$.  Thus, $f = 0$. 
It follows that this case yields the irreducible component 
\[
\L_2 = \V(\, M_{13},\,
M_{14}, \ M_{34}, \ M_{12}^3 - M_{12} M_{23}^2 -i M_{23} M_{24}^2\, )
\]
of $\L'(\g)$.

\bigskip

\noindent
\textbf{Case (III):}
$M_{14}M_{24} \neq 0$ and $M_{23}=0$.\\
This case is computed by applying $\psi_1$ to case (II), giving 
\[
\L_3= \V(\, M_{12}, \ M_{13}, \ M_{23}, \ M_{34}^3 - M_{14}^2 M_{34} 
+ i M_{14} M_{24}^2\, ).
\]

\bigskip

\noindent
\textbf{Case (IV):}
$M_{23} \neq 0$ and $M_{14}=0=M_{24}$.\\
If, additionally, $M_{12} \neq 0$, then $M_{13} = 0$ and $M_{i4} = 0$
for all $i = 1,\,2,\,3$. It follows that $M_{12}^2 = M_{23}^2$, and
so these assumptions yield a subvariety of $\L_2$. Hence, we may assume
that $M_{12} = 0$. It follows that this case yields the irreducible 
component 
\[
\L_4 = \V(\,M_{12}, \ M_{14}, \ M_{24}, \ M_{23}^2 M_{34} + i \g M_{13}^2
M_{23} - M_{34}^3\,)
\]
of $\L'(\g)$, so $\L_4$ is $\psi_2$ applied to $\L_2$.

\bigskip

\noindent
\textbf{Case (V):}
$M_{14} \neq 0$ and $M_{23}=0=M_{24}$.\\
This case is computed by applying $\psi_1$ to case (IV), giving the 
irreducible component 
\[
\L_5
= \V(\,M_{23}, \ M_{24}, \ M_{34}, \ M_{12} M_{14}^2  - i \g M_{13}^2 M_{14} 
- M_{12}^3\,)
\]
of $\L'(\g)$, which is also $\psi_2$ applied to $\L_3$.

\bigskip

\noindent
\textbf{Case (VI):}
$M_{24} \neq 0$ and $M_{14}=0=M_{23}$.\\
Using $M_{14} = 0 = M_{23}$, a computation of a Gr\"obner basis yields 
several polynomials, one of which is $M_{12} M_{34}- M_{13} M_{24}$ whereas 
the others are multiples of $M_{12}^2 + M_{34}^2$. In particular, two
of those polynomials are: $M_{12} M_{24} (M_{12}^2 + M_{34}^2)$ and 
$M_{34}^2 (M_{12}^2 + M_{34}^2)$.  It follows that $M_{12}^2 + M_{34}^2
= 0$, so that this case yields the component $\L_6 = 
\L_{6a} \cup \L_{6b}$ of $\L'(\g)$, where 
\[
\begin{array}{c}
\L_{6a} = 
\V(\, M_{14}, \ M_{23}, \ M_{12} M_{34}- M_{13} M_{24}, \ M_{12} + i M_{34}
\,), \\[2mm]
\L_{6b} = 
\V(\, M_{14}, \ M_{23}, \ M_{12} M_{34}- M_{13} M_{24}, \ M_{12} - i M_{34}
\,),
\end{array}
\]
and each of $\L_{6a}$ and $\L_{6b}$ is a nonsingular conic, since 
using $M_{12} \pm i M_{34}$ to substitute for $M_{12}$ in 
$M_{12} M_{34}- M_{13} M_{24}$ yields a rank-3 quadratic form in each case.
Moreover, $\L_{6b}$ is $\psi_1$ applied to $\L_{6a}$.

\bigskip

Having completed this analysis, we can see that the point 
$\V(\, M_{12},\, M_{14},\, M_{23},\, M_{24},\, M_{34}\,)$, that was found 
earlier, is contained in $\L_4\cap \L_5 \cap \L_6$. 
We summarize the above work in the next result.

\begin{thm} \label{L'g}
Let $\L'(\g)$ denote the reduced variety of the line scheme 
$\L(\g)$ of $\A(\g)$.
If $\g^2 \neq 16$, then $\L'(\g)$ is the union, in~$\PP^5$, of the 
following seven irreducible components:
\begin{enumerate}
\item[\textup{(I)}] $\L_1 = \V(\, M_{13}, \ M_{24}, \ 
              M_{14} M_{23}+M_{12} M_{34}, \ 
        M_{12}^2 + M_{34}^2 + \g M_{14} M_{23} - M_{14}^2 - M_{23}^2 \,),$ 
	      which is a nonplanar elliptic curve in a $\PP^3$. 
\item[\textup{(II)}] $\L_2  = \V(\, M_{13},\, M_{14}, \ M_{34}, \ 
              M_{12}^3 - M_{12} M_{23}^2 -i M_{23} M_{24}^2\, ),$
              which is a planar elliptic curve.
\item[\textup{(III)}] $\L_3= \V(\, M_{12}, \ M_{13}, \ M_{23}, \ 
              M_{34}^3 - M_{14}^2 M_{34} + i M_{14} M_{24}^2\, ),$ 
              which is a planar elliptic curve.
\item[\textup{(IV)}] $\L_4 = \V(\,M_{12}, \ M_{14}, \ M_{24}, \ 
              M_{23}^2 M_{34} + i \g M_{13}^2 M_{23} - M_{34}^3\,),$ 
              which is a planar elliptic curve.
\item[\textup{(V)}] $\L_5 = \V(\,M_{23}, \ M_{24}, \ M_{34}, \ 
              M_{12} M_{14}^2  - i \g M_{13}^2 M_{14} - M_{12}^3\,),$ 
              which is a planar elliptic curve.
\item[\textup{(VIa)}] $\L_{6a} = \V(\, M_{14}, \ M_{23}, \ 
              M_{12} M_{34}- M_{13} M_{24}, \ M_{12} + i M_{34} \,),$ 
	      which is a nonsingular conic.
\item[\textup{(VIb)}] $\L_{6b} = \V(\, M_{14}, \ M_{23}, \ 
              M_{12} M_{34}- M_{13} M_{24}, \ M_{12} - i M_{34}\,),$ 
	      which is a nonsingular conic.
\end{enumerate}
If $\g = 4$, then $\L'(\g)$ is the union, in~$\PP^5$, of eight irreducible
components, six of which are $\L_2$, $\L_3$, $\L_4$, 
$\L_5$, $\L_{6a}$, $\L_{6b}$ {\em (}\!as above{\em )}
and two of which are 
\[
\begin{array}{c}
\L_{1a} = \V(\, M_{13}, \ M_{24}, \ M_{14} M_{23}+M_{12} M_{34},  \ 
             M_{12} + M_{14} - M_{23} - M_{34}\,),\\[2mm]
\L_{1b} = \V(\, M_{13}, \ M_{24}, \ M_{14} M_{23}+M_{12} M_{34},  \ 
             M_{12} - M_{14} + M_{23} - M_{34}\,),
\end{array}
\]
which are nonsingular conics.
\end{thm}
\begin{pf}
The polynomials were found in the preceding work, as was the geometric 
description for $\L_{1a}$, $\L_{1b}$, $\L_{6a}$ and $\L_{6b}$, so here we 
discuss only the geometric description of the other components.

(I) Write $q_1 = M_{14} M_{23}+M_{12} M_{34}$ and $q_2 = 
M_{12}^2 + M_{34}^2 + \g M_{14} M_{23} - M_{14}^2 - M_{23}^2$ 
viewed in $\kk[M_{12}, \, M_{14}, \, M_{23}, \, M_{34}]$.  Since 
\[
q_2 = M_{12}^2 - (\g/2) M_{12} M_{34} + M_{34}^2 -
\left( M_{14}^2 - (\g/2) M_{14} M_{23} + M_{23}^2 \right)
\]
modulo $q_1$, and since 
char$(\kk) \neq 2$, we may take the Jacobian matrix of this system of two 
polynomials to be the $2 \times 4$ matrix\\[-2mm]
\[
\footnotesize
\begin{bmatrix}
M_{34} & M_{23} & M_{14} & M_{12} \\[3mm]
2 M_{12} - (\g/2) M_{34} & 
-(2 M_{14} - (\g/2) M_{23}) & 
-(2 M_{23} - (\g/2) M_{14}) & 
2 M_{34} - (\g/2) M_{12} & 
\end{bmatrix}.
\]
\quad\\[-1mm]
Assuming that all the $2 \times 2$ minors are zero, we find that 
$M_{34}^2 = M_{12}^2$ (from columns one and four) and 
$M_{23}^2 = M_{14}^2$ (from columns two and three). Substituting these
relations into the minor obtained from the last two columns yields
either $(\g \pm 4)M_{12} M_{14} = 0$ or $\g M_{12} M_{14} = 0$, so 
$ M_{12} M_{14}= 0$ (since $\g (\g^2-16)\neq 0$). 
Substitution into $q_1$ implies that there is no solution, and so the 
Jacobian matrix has rank two at all points of $\V(q_1,\, q_2)$. 
It follows that $\V(q_1,\, q_2)$, viewed as a
subvariety of $\PP^3 = \V(M_{13},\, M_{24})$, is reduced, and so
$\L_1$ is reduced. Following the method of the proof of
\cite[Proposition~2.5]{Smith.Staff}, if $\V(q_1,\, q_2)$ is not irreducible,
then there exists a point in the intersection of two of its irreducible
components, and so the Jacobian matrix has rank at most one at that point, 
which is a contradiction. Hence, $\V(q_1,\, q_2)$ is irreducible, and thus 
nonsingular since it is reduced. Moreover, its genus is $4-2-2+1 = 1$. 
It follows that $\V(q_1,\, q_2)$ is an elliptic curve, and the same is true 
of $\L_1$. 

(II) Viewing $h = M_{12}^3 - M_{12} M_{23}^2 -i M_{23} M_{24}^2$ as a 
polynomial in $\kk[M_{12}, \, M_{23}, \, M_{24}]$, the Jacobian matrix of 
$h$ is a $1 \times 3$ matrix that has rank one at all points of $\V(h)$ 
(since char$(\kk) \neq 2$), so $\V(h)$ is nonsingular in $\PP^2 = 
\V(M_{13},\, M_{14}, \ M_{34})$.

(III), (IV), (V) These cases follow from (II) by applying $\psi_1$ or
$\psi_2$ as appropriate.
\end{pf}

\bigskip
%\newpage

\subsection{Description of the Line Scheme}\hfill

In this subsection, we prove that the line scheme $\L(\g)$ of $\A(\g)$ is
reduced and so is given by $\L'(\g)$ described in Theorem~\ref{L'g}.

\begin{lemma}\label{dim1}
For all $\g \in \kk^\times$, the irreducible components of $\L(\g)$ have 
dimension one; in particular, $\L(\g)$ has no embedded points.
\end{lemma}
\begin{pf}
By \cite{CV}, $\A(\g)$ is a regular noetherian domain that is 
Auslander-regular and satisfies the Cohen-Macaulay property and has
Hilbert series the same as that of the polymomial ring on four variables.
Hence, by \cite[Remark~2.10]{SV1}, we may apply \cite[Corollary~2.6]{SV1} 
to $\A(\g)$, which gives us that the irreducible components of $\L(\g)$ 
have dimension at least one. However, by Theorem~\ref{L'g}, they have
dimension at most one, so equality follows. 
Let $X_1$ denote the 11-dimensional subscheme of $\PP(V \otimes V)$
consisting of the elements of rank at most two, and, for all $\g \in
\kk^\times$, let $X_2$ denote the 5-dimensional linear subscheme of 
$\PP(V \otimes V)$ given by the span of the defining relations of 
$\A(\g)$. By \cite[Lemma~2.5]{SV1}, $\L(\g)\cong X_1 \cap X_2$ for
all $\g \in \kk^\times$.  Since $X_i$ is a Cohen-Macaulay scheme for 
$i = 1, \, 2$, and since $\dim(X_1 \cap X_2) = 1$,  the proof of 
\cite[Theorem~4.3]{SV1} (together with Macaulay's Unmixedness Theorem) 
rules out the possibility of embedded components.
\end{pf}

\begin{thm}\label{Lg}
For all $\g \in \kk^\times$, the line scheme~$\L(\g)$ is a reduced scheme 
of degree twenty.
\end{thm}
\begin{pf}
Let $X_1$ and $X_2$ be as in the proof of Lemma~\ref{dim1}, and let
$X = X_1 \cap X_2$. Since $\deg(X_1) = 20$ by \cite[Example~19.10]{H}, 
B\'ezout's Theorem for Cohen-Macaulay schemes (\cite[Theorem~III-78]{EH})
implies that $\deg(X) = 20$.
However, since $\L(\g) \cong X$ by \cite[Lemma~2.5]{SV1}, 
the reduced scheme $X'$ of $X$ is isomorphic to $\L'(\g)$.
Since the degrees of the irreducible components of $\L'(\g)$
in Theorem~\ref{L'g} are as small as possible, 
$\deg(X') \geq 4 + 12 + 4 = 20$; that is,
$20 = \deg(X) \geq \deg(X') \geq 20$,  giving $\deg(X) = \deg(X')$. 
As $X$ has no embedded points by Lemma~\ref{dim1}, it follows that 
$X = X'$, so $X$ is a reduced scheme. Thus, $\L(\g)$ is reduced and has 
degree twenty since $\deg(\L'(\g)) = 20$.
\end{pf}

The intersection points of the irreducible components of $\L(\g)$ are 
straightforward to compute and are listed in \cite{C}.

\bigskip
%\newpage

%%%%%%%%%%%%%%%%%%%%%%%%%%%%%%%%%%%%%%%%%%%%%%%%%%%%%%%%%%%%%%%%%%%%%
% END  SECTION  3
% START  SECTION  4
%%%%%%%%%%%%%%%%%%%%%%%%%%%%%%%%%%%%%%%%%%%%%%%%%%%%%%%%%%%%%%%%%%%%%

\section{The Lines in $\PP^3$ Parametrized by the Line Scheme}\label{sec4}

In this section, we describe the lines in $\PP(V^{\text{*}})$ that are
parametrized by the line scheme $\L(\g)$ of $\A(\g)$.   We also
describe, in Theorem~\ref{6lines}, the lines that pass through any given 
point of the point scheme; in particular, if $p$ is one of the generic 
points of the point scheme (that is, $p \in \Z_\gamma$), then there are 
exactly six distinct lines of the line scheme that pass through~$p$.
Since we will use results from Section~\ref{sec3}, we continue to assume 
that char$(\kk) = 0$.

\medskip

\subsection{The Lines in $\PP^3$}\hfill

In this subsection, we find the lines in $\PP(V^{\text{*}})$ that are
parametrized by the line scheme.
We first recall how the Pl\"ucker coordinates $M_{12},\ldots, M_{34}$ relate
to lines in $\PP^3$; details may be found in \cite[\S8.6]{CLO}.
Any line $\ell$ in $\PP^3$ is uniquely determined by any two distinct points
$a = (a_1,\ldots ,a_4) \in \ell$ and $b = (b_1,\ldots , b_4) \in \ell$, and 
may be represented by a $2 \times 4$ matrix
\[
\begin{bmatrix}
a_1 & a_2 & a_3 & a_4 \\[2mm]
b_1 & b_2 & b_3 & b_4 
\end{bmatrix}
\]
that has rank two; in particular, the points on $\ell$ are represented
in homogeneous coordinates by linear combinations of the rows of this matrix.
In general, there are infinitely many such matrices that may 
be associated to any line $\ell$ in $\PP^3$, and they are all related to 
each other by applying row operations.

The Pl\"ucker coordinate $M_{ij}$ is evaluated on this matrix as the minor 
$a_i b_j - a_j b_i$ for all $i \neq j$, and the Pl\"ucker polynomial 
$P = M_{12}M_{34}-M_{13}M_{24}+M_{14}M_{23}$, given in
Section~\ref{method}, vanishes on this matrix. Moreover, $\V(P)$ is the
subscheme of $\PP^5$ that parametrizes all lines in $\PP^3$.

Since $\dim(V) =4$, we identify $\PP(V^{\text{*}})$ with $\PP^3$.
By Theorem~\ref{Lg}, $\L(\g)$ is given by Theorem~\ref{L'g}. We
continue to use the notation $e_j$ introduced in Section~\ref{sec2}.

\bigskip

\noindent \textbf{(I)}
 \ In this case, $\g^2 \neq 16$ and the component is $\L_1$, 
which is a nonplanar elliptic curve in a $\PP^3$ (contained in $\PP^5$), 
where  
\[
\L_1 = \V(\, M_{13}, \ M_{24}, \ M_{14} M_{23}+M_{12} M_{34}, \ 
M_{12}^2 + M_{34}^2 + \g M_{14} M_{23} - M_{14}^2 - M_{23}^2
\, ).
\]
It follows that any line $\ell$ in $\PP(V^{\text{*}})$ given by $\L_1$ is 
represented by a $2 \times 4$ matrix of the form:\\[-1mm]
\[
\begin{bmatrix}
a_1 & 0 & a_3 & 0 \\[2mm]
0 & b_2 & 0 & b_4 
\end{bmatrix},
\tag{$*$}
\]
\quad\\
where $a_j,\,b_j \in \kk$ for all $j$ and 
$a_1^2 b_2^2 + a_3^2 b_4^2 -\g a_1 b_2 a_3 b_4 - a_1^2 b_4^2 - b_2^2
a_3^2 = 0$. In particular, if $p \in \ell$, then 
$p = (\lambda_1 a_1,\, \lambda_2 b_2,\, \lambda_1 a_3,\, \lambda_2 b_4)$,
for some $(\lambda_1 ,\, \lambda_2 )\in \PP^1$, such that 
$a_1^2 b_2^2 + a_3^2 b_4^2 -\g a_1 b_2 a_3 b_4 - a_1^2 b_4^2 - b_2^2
a_3^2 = 0$.  It is easily verified that $p$ lies on the quartic surface 
\[ \V(\, 
x_1^2 x_2^2 + x_3^2 x_4^2 -\g x_1 x_2 x_3 x_4 - x_1^2 x_4^2 - x_2^2 x_3^2 
\,)
\]
in $\PP(V^{\text{*}})$ for all $(\lambda_1 ,\, \lambda_2 )\in \PP^1$.  
Hence, the lines parametrized by $\L_1$ all lie on this quartic surface
in $\PP(V^{\text{*}})$ and are given by: 
\[
\V(x_3,\ x_2 \pm x_4), \quad 
\V(x_4,\ x_1 \pm x_3), \quad \text{and} \quad
\V(x_1 - \alpha x_3, \ x_2 -\beta x_4)
\]
for all $\alpha,\ \beta \in \kk$ such that 
$(\alpha^2-1)(\beta^2 -1) = \g \alpha \beta$. The case $\g = 4$ is
discussed below.

\bigskip

\noindent \textbf{(II)}
 \ In this case, the component is $\L_2$, which is a planar 
elliptic curve, where 
\[
\L_2  = \V(\, M_{13},\, M_{14}, \ M_{34}, \ 
              M_{12}^3 - M_{12} M_{23}^2 -i M_{23} M_{24}^2\, ), 
\]
so any line in $\PP(V^{\text{*}})$ given by $\L_2$ is represented by a 
$2 \times 4$ matrix of the form:\\[-1mm]
\[
\begin{bmatrix}
a_1 & 0 & a_3 & a_4 \\[2mm]
0 & 1 & 0 & 0 
\end{bmatrix},
\]
\quad\\[-1mm]
such that $a_1^3 - a_1 a_3^2 + i a_3 a_4^2 = 0$. It follows that $\L_2$
parametrizes those lines in $\PP(V^{\text{*}})$ that pass through $e_2$
and meet the planar curve 
$\V(x_2, \, x_1^3 - x_1 x_3^2 + i x_3 x_4^2 )$; this planar curve is a
(nonsingular) elliptic curve since char$(\kk) = 0$.

\bigskip

\noindent \textbf{(III)}
 \ In this case, the component is $\L_3$, which may be
obtained as $\psi_1$  applied to $\L_2$. Hence, $\L_3$ parametrizes those 
lines in $\PP(V^{\text{*}})$ that pass through $e_4$ and 
meet the planar elliptic curve 
$\V(x_4, \, x_3^3 - x_1^2 x_3 + i x_1 x_2^2 )$. 

\bigskip

\noindent \textbf{(IV)}
 \ In this case, the component is $\L_4$, which may be
obtained as $\psi_2$  applied to $\L_2$. Hence, $\L_4$ parametrizes those 
lines in $\PP(V^{\text{*}})$ that pass through $e_3$ and 
meet the planar elliptic curve 
$\V(x_3, \, x_4^3 - x_2^2 x_4 + i \g x_1^2 x_2)$. 

\bigskip

\noindent \textbf{(V)}
 \ In this case, the component is $\L_5$, which may be
obtained as $\psi_1$  applied to $\L_4$. Hence, $\L_5$ parametrizes those 
lines in $\PP(V^{\text{*}})$ that pass through $e_1$ and 
meet the planar elliptic curve 
$\V(x_1, \, x_2^3 - x_2 x_4^2 + i \g x_3^2 x_4)$. 

\bigskip

\noindent \textbf{(VI)}
 \ In this case, the component is 
$\L_6 = \L_{6a} \cup \L_{6b}$, where 
\[
\begin{array}{c}
\L_{6a} = \V(\, M_{14}, \ M_{23}, \ 
              M_{12} M_{34}- M_{13} M_{24}, \ M_{12} + i M_{34} \,),\\[2mm]
\L_{6b} = \V(\, M_{14}, \ M_{23}, \ 
              M_{12} M_{34}- M_{13} M_{24}, \ M_{12} - i M_{34}\,), 
\end{array}
\]
which are nonsingular conics. Following the argument from case (I), any line 
in $\PP(V^{\text{*}})$ given by $\L_{6a}$ is represented by a $2 \times 4$ 
matrix of the form:\\[-1mm]
\[
\begin{bmatrix}
a_1 & a_2 & a_3 & a_4 \\[2mm]
\alpha a_1 & \beta a_2 & \beta a_3 & \alpha a_4
\end{bmatrix},
\]
\quad\\[-1mm]
such that $\alpha,\, \beta,\, a_j \in \kk$ for all $j$, 
$a_1 a_2 = i a_3 a_4$
and $\alpha \neq \beta$. A calculation similar to that used in~(I) 
verifies that every point of the
line lies on the quadric $\V(x_1 x_2 - i x_3 x_4)$.  
It follows that $\L_{6a}$ parametrizes one of the 
rulings of the nonsingular quadric $\V(x_1 x_2 - i x_3 x_4)$; namely, the 
ruling that consists of the lines 
$\V(\delta x_1 - \epsilon x_4,\, \delta x_3 + i \epsilon x_2 )$ for all 
$(\delta, \, \epsilon ) \in \PP^1$. Since $\L_{6b}$ may be obtained by
applying $\psi_1$ to $\L_{6a}$, we find $\L_{6b}$ parametrizes one of the 
rulings of the nonsingular quadric $\V(x_3 x_4 - i x_1 x_2)$; namely, the 
ruling that consists of the lines 
$\V(\delta x_3 - \epsilon x_2,\, \delta x_1 + i \epsilon x_4 )$ for all
$(\delta, \, \epsilon ) \in \PP^1$.

\bigskip
\medskip

\noindent \textbf{(Ia)} and \textbf{(Ib)} 
 \ In this case, $\g = 4$ and the component is 
$\L_1 = \L_{1a} \cup \L_{1b}$, where 
\[
\begin{array}{c}
\L_{1a} = \V(\, M_{13}, \ M_{24}, \ M_{14} M_{23}+M_{12} M_{34},  \ 
             M_{12} + M_{14} - M_{23} - M_{34}\,),\\[2mm]
\L_{1b} = \V(\, M_{13}, \ M_{24}, \ M_{14} M_{23}+M_{12} M_{34},  \ 
             M_{12} - M_{14} + M_{23} - M_{34}\,),
\end{array}
\]
which are nonsingular conics. Following the argument from case (I), any line 
in $\PP(V^{\text{*}})$ given by $\L_{1a}$ is represented by a $2 \times 4$ 
matrix of the form \thetag{$*$} such that $a_1 b_2 + a_1 b_4 + b_2 a_3 
= a_3 b_4$. A calculation similar to that used in~(I) verifies that every 
point of the line lies on the nonsingular quadric 
\[  Q_a = \V(\, x_1 x_2 + x_1 x_4 + x_2 x_3 - x_3 x_4 \,) \]
in $\PP(V^{\text{*}})$. Hence, the lines parametrized by $\L_{1a}$ all
lie on $Q_a$ and are: 
\[
\V(x_3,\ x_2 + x_4) \quad \text{and} \quad
\V(x_1 - \alpha x_3, \ (\alpha + 1) x_2 + (\alpha - 1) x_4)
\]
for all $\alpha \in \kk$, which yields one of the rulings on the quadric
$Q_a$. Applying $\psi_1$ to these lines, it follows that the lines 
parametrized by $\L_{1b}$ are: 
\[
\V(x_1,\ x_2 + x_4) \quad \text{and} \quad
\V(x_3 - \alpha x_1, \ (\alpha - 1) x_2 + (\alpha + 1) x_4)
\]
for all $\alpha \in \kk$, which yields one of the rulings on the
nonsingular quadric
\[ Q_b = \V( \, x_3 x_4 + x_2 x_3 + x_1 x_4 - x_1 x_2 \,).\] 

\bigskip

\subsection{The Lines of the Line Scheme That Contain Points of the Point 
Scheme}\hfill

In this subsection, we compute how many lines in $\PP(V^{\text{*}})$ that 
are parametrized by~$\L(\g)$ contain a given point of~$\p(\g)$. 
By \cite[Remark~3.2]{SV1}, if the number of lines is finite,
then it is six, counting multiplicity; hence, the generic case is considered
to be six distinct lines. 
The reader should note that a result similar to Theorem~\ref{6lines} is 
given in \cite[Theorem~IV.2.5]{G} for the algebra $\A(1)$, but that result 
is false as stated (perhaps as a consequence of the sign errors in the
third relation of (3) on Page 797 of \cite{ST}). 

\begin{thm}\label{6lines}
Suppose $\g \in \kk^\times$, and let $\Z_\g$ be as in 
Proposition~\ref{ptscheme}.
\begin{enumerate}
\itemsep1mm
\item[\textup{(a)}]
     For any $j \in \{ 1, \ldots, 4\}$, $e_j$ lies on infinitely many lines 
     that are parametrized by $\L(\g)$.
\item[\textup{(b)}]
     Each point of $\Z_\g$ lies on exactly six distinct lines of
     those parametrized by $\L(\g)$.
\end{enumerate}
\end{thm}
\begin{pf}
Since (a) follows from (II)-(V) in Section~4.1, we focus on (b).  Let $p =
(1,\, \alpha_2,\, \alpha_3,\, \alpha_4) \in \Z_\g$. It follows that
$\alpha_j \neq 0$ for all $j$.  Suppose that $\g^2 \neq 16$.

Let $\alpha = 1/\alpha_3$ and $\beta = \alpha_2 /\alpha_4$, so 
$(\alpha^2 -1)(\beta^2 -1) = \g \alpha \beta$, by \ref{app1}.\ref{pp15} 
in Section~\ref{app1}. Hence, $p \in \V(x_1 - \alpha x_3, \ x_2 -\beta x_4)$,
which is a line that corresponds to an element of $\L_1$.  Clearly, no other
line given by $\L_1$ contains $p$.

Let $r_2 = (1,\, 0,\, \alpha_3,\, \alpha_4)$ and let $\ell_2$ denote the
line through $e_2$ and $r_2$.  By \ref{app1}.\ref{pp9}, we have  
$1 - \alpha_3^2 + i \alpha_3 \alpha_4^2 = 0$, so $r_2 \in
\V(x_2, \, x_1^3 - x_1 x_3^2 + i x_3 x_4^2 )$.  Thus, $\ell_2$ 
corresponds to an element of $\L_2$, and $p \in \ell_2$.  Conversely,
let $r_2' = (b_1,\, 0,\, b_3,\, b_4) \in \V(x_2, \, x_1^3 - x_1 x_3^2 +
i x_3 x_4^2 )$. If $p$ lies on the line through $r_2'$ and $e_2$, then
there exist $(\lambda_1, \, \lambda_2) \in \PP^1$ such that $p = 
(\lambda_1 b_1,\, \lambda_2 ,\, \lambda_1 b_3,\, \lambda_1 b_4)$. Thus,
$\lambda_1 b_1 \neq 0$ and $\alpha_i = b_i/b_1$ for $i = 3,\, 4$.
Hence, $r_2' = (b_1,\, 0,\, b_1 \alpha_3,\, b_1 \alpha_4) = 
(1,\, 0,\, \alpha_3,\, \alpha_4) = r_2$. It follows that no other line given
by $\L_2$ contains $p$.

Let $r_4 = (1,\, \alpha_2,\, \alpha_3,\, 0)$ and let $\ell_4$ denote the
line through $e_4$ and $r_4$.  By \ref{app1}.\ref{pp2}, we have  
$\alpha_3^3 - \alpha_3 + i \alpha_2^2 = 0$, so $r_4 \in
\V(x_4, \, x_3^3 - x_1^2 x_3 + i x_1 x_2^2 )$.  Thus, $\ell_4$ 
corresponds to an element of $\L_3$, and $p \in \ell_4$. An argument 
similar to that of $\L_2$ proves that no other line given by $\L_3$ 
contains~$p$.  

Let $r_3 = (1,\, \alpha_2,\, 0,\, \alpha_4 )$ and let $\ell_3$ denote the
line through $e_3$ and $r_3$.  By \ref{app1}.\ref{pp5}, we have  
$\alpha_4^3- \alpha_2^2 \alpha_4 + i \g \alpha_2$ = 0, so $r_3 \in
\V(x_3, \, x_4^3- x_2^2 x_4 + i \g x_1^2 x_2)$.  Thus, $\ell_3$ 
corresponds to an element of $\L_4$, and $p \in \ell_3$. 
An argument similar to that of $\L_2$ proves that no other line given by 
$\L_4$ contains~$p$.  

Let $r_1 = (0,\, \alpha_2,\, \alpha_3,\, \alpha_4 )$ and let $\ell_4$ denote
the line through $e_1$ and $r_1$.  By \ref{app1}.\ref{pp8}, we have  
$\alpha_2^3- \alpha_2 \alpha_4^2 + i \g \alpha_3^2 \alpha_4 = 0$,
so $r_1 \in \V(x_1, \, x_2^3- x_2 x_4^2 + i \g x_3^2 x_4 )$.  Thus, $\ell_4$ 
corresponds to an element of $\L_5$, and $p \in \ell_4$. 
An argument similar to that of $\L_2$ proves that no other line given by 
$\L_5$ contains~$p$.  

By \ref{app1}.\ref{pp1}, we have $\alpha_2 = \pm i \alpha_3 \alpha_4$, 
so either $p \in \V(x_1 x_2 - i x_3 x_4)$ or 
$p \in \V(i x_1 x_2 - x_3 x_4)$ (but not both, since $\alpha_3 \alpha_4
\neq 0$). In the first case, $p \in 
\V(\alpha_4 x_1 - x_4,\, \alpha_4 x_3 + i x_2)$ and, in the second,
$p \in \V( \alpha_4 x_1 - x_4,\, i \alpha_4 x_3 + x_2)$.  These lines
correspond to elements of $\L_{6a}$ and $\L_{6b}$ respectively. Since
each quadric has only two rulings, and since each irreducible component
of $\L_6$ parametrizes 
only one of the rulings in each case, no other line given by $\L_6$ 
contains~$p$.  

If, instead, $\g = 4$, the only adjustment to the above reasoning is in 
the case of the lines parametrized by $\L_1$.  Since $\g = 4$, 
the polynomial \ref{app1}.\ref{pp15} factors, so
\[
(\alpha_2 + \alpha_4 + \alpha_2 \alpha_3 - \alpha_3 \alpha_4)
(\alpha_2 - \alpha_4 - \alpha_2 \alpha_3 - \alpha_3 \alpha_4) = 0,
\tag{$\dag$}
\]
that is, 
\[
\big((1+ \alpha_3)\alpha_2  + (1 - \alpha_3 )\alpha_4 \big)
\big((1- \alpha_3)\alpha_2  - (1 + \alpha_3 )\alpha_4 \big) = 0,
\]
which provides exactly two lines (of those parametrized by $\L_1$) that 
could contain~$p$. These lines are 
\[
\V(x_1 - (1/\alpha_3) x_3, \, ((1/\alpha_3)+1)x_2 + (1/\alpha_3)-1)x_4)
\]
and 
\[
\V(x_3 - \alpha_3 x_1, \, (\alpha_3-1)x_2 + (\alpha_3+1)x_4),
\]
which correspond to elements of $\L_{1a}$ and $\L_{1b}$ respectively.
If the first factor of \thetag{$\dag$} is zero, then $p$ belongs to the
first line, whereas if the second factor of \thetag{$\dag$} is zero, then
$p$ belongs to the second line. If both factors of \thetag{$\dag$} 
are zero, then 
$\alpha_2 = \alpha_3 \alpha_4$, which forces $\alpha_3 \alpha_4 = 0$, by 
\ref{app1}.\ref{pp1}, and this contradicts $p \in \Z_\g$. It follows that
$p$ belongs to exactly one line of those parametrized by $\L_1$.

For all $\gamma \in \kk^\times$, it is a straightforward calculation to 
show that the six lines found above are distinct. 
\end{pf}

\bigskip

Considering Theorems~\ref{L'g}, \ref{Lg} and \ref{6lines} in the case 
where $\g^2 \neq 16$, we arrive at the following conjecture.

\begin{conj}\label{conj}
The line scheme of the most generic quadratic quantum~$\PP^3$ is isomorphic
to the union of two spatial (irreducible and nonsingular) elliptic curves 
and four planar (irreducible and nonsingular) elliptic curves. 
(Here, {\em spatial elliptic curve} means a nonplanar elliptic curve that
is contained in a subscheme of~$\PP^5$ that is isomorphic to~$\PP^3$.)
\end{conj}

This conjecture is motivated by the idea that the ``generic'' points of
the point scheme should have exactly six distinct lines of the line scheme 
passing 
through each of them, with each line coming from exactly one component of 
the line scheme. Moreover, if the component $\L_6$ of the line scheme
$\L(\g)$ of $\A(\g)$ had not split into two smaller components,
then it would likely have been a spatial elliptic curve.

%%%%%%%%%%%%%%%%%%%%%%%%%%%%%%%%%%%%%%%%%%%%%%%%%%%%%%%%%%%%%%%%%%%%%
% END  SECTION  4
% START  APPENDIX (SEC 5)
%%%%%%%%%%%%%%%%%%%%%%%%%%%%%%%%%%%%%%%%%%%%%%%%%%%%%%%%%%%%%%%%%%%%%

\bigskip
\bigskip
%\newpage

\section{Appendix}\label{app}

In this section, we list the polynomials that define $\p(\g)$ and 
$\L(\g)$.

\medskip

\subsection{Polynomials Defining the Point Scheme}\label{app1}\hfill

The following are the polynomials that define the point scheme viewed as 
$\p(\g)\subset \PP(V^{\text{*}})$ of 
$\A(\g)$ that are given by the fifteen $4 \times 4$ minors of the 
matrix~$M$ in Section~\ref{sec2};
they are used in Section~\ref{sec2} and in the proof of 
Theorem~\ref{6lines}:
\begin{enumerate}[\qquad\thesubsection.1.\quad]
\itemsep1mm
\item\label{pp1}
$x_1^2 x_2^2 + x_3^2 x_4^2$, 
\item\label{pp2}
$x_1 \left( x_3^3- x_1^2 x_3 + i x_1 x_2^2 \right)$, 
\item
$x_2 \left( x_3^3- x_1^2 x_3 + i x_1 x_2^2 \right)$, 
\item
$x_4 \left( x_3^3- x_1^2 x_3 + i x_1 x_2^2 \right)$, 
\item\label{pp5} 
$x_1 \left( x_4^3- x_2^2 x_4 + i \g x_1^2 x_2 \right)$, 
\item
$x_2 \left( x_4^3- x_2^2 x_4 + i \g x_1^2 x_2 \right)$, 
\item
$x_3 \left( x_4^3- x_2^2 x_4 + i \g x_1^2 x_2 \right)$, 
\item\label{pp8}
$x_1 \left( x_2^3- x_2 x_4^2 + i \g x_3^2 x_4 \right)$, 
\item\label{pp9}
$x_2 \left( x_1^3- x_1 x_3^2 + i x_3 x_4^2 \right)$, 
\item\label{pp10}
$i \g x_1^2 x_3^2 - x_1^2 x_2 x_4 - x_2 x_3^2 x_4$, 
\item
$i x_2^2 x_4^2 - x_1 x_2^2 x_3 - x_1 x_3 x_4^2$, 
\item
$x_1^3 x_4 + \g x_1^2 x_2 x_3 - x_1 x_3^2 x_4 + i x_2^2 x_3
x_4$,
\item
$x_2^3 x_3 + \g x_1 x_2^2 x_4 - x_2 x_3 x_4^2 + i \g x_1^2
x_3 x_4$,
\item
$i \g x_1^3 x_3 + \g x_1^2 x_2^2 - 2 x_1 x_2 x_3 x_4 + i x_2^3
x_4$,
\item\label{pp15}
$x_1^2 x_2^2 - x_2^2 x_3^2 - \g x_1 x_2 x_3 x_4 - x_1^2 x_4^2 +
        x_3^2 x_4^2$,
\end{enumerate}
\quad\\[-3mm]
where $i^2 = -1$ and $\g \in \kk^\times$.

\bigskip

\subsection{Polynomials Defining the Line Scheme}\label{app2}\hfill

The following are the forty-six polynomials in the $M_{ij}$ coordinates from 
Section~\ref{sec3} that define the line scheme $\L(\g)$ of $\A(\g)$:
\begin{enumerate}[\qquad\thesubsection.1.\quad]
\itemsep1mm
\setcounter{enumi}{-1}
\item
$P = M_{12} M_{34} - M_{13} M_{24} + M_{14} M_{23}$,

\item 
$ 2 M_{13} M_{14} M_{23} M_{24} $,

\item 
$ M_{12} (\g M_{13} M_{14} M_{23}+i M_{12} M_{14} M_{24}+i M_{23} M_{24} M_{34}) $,

\item 
$ M_{12} (\g M_{13} M_{14} M_{23}-i M_{12} M_{14} M_{24}-i M_{23} M_{24} M_{34}) $,

\item 
$ M_{13} (\g M_{13} M_{14} M_{23}+i M_{12} M_{14} M_{24}+i M_{23} M_{24} M_{34}) $,

\item 
$ M_{13} (\g M_{13} M_{14} M_{23}-i M_{12} M_{14} M_{24}-i M_{23} M_{24} M_{34}) $,

\item 
$M_{13} (\g M_{13} M_{14} M_{23}+i M_{12} M_{14} M_{24}-i M_{23} M_{24} M_{34}) $,

\item 
$ M_{14} (\g M_{13} M_{14} M_{23}+i M_{12} M_{14} M_{24}+i M_{23} M_{24} M_{34}) $,

\item 
$ M_{23} (\g M_{13} M_{14} M_{23}+i M_{12} M_{14} M_{24}+
     i M_{23} M_{24} M_{34}) $,

\item 
$M_{23} (\g M_{13} M_{14} M_{23}-i M_{12} M_{14} M_{24}-i M_{23} M_{24}
 M_{34}) $,

\item 
$ M_{24} (\g M_{13} M_{14} M_{23}+i M_{12} M_{14} M_{24}+i M_{23}
 M_{24} M_{34}) $,

\item 
$ M_{34} (\g M_{13} M_{14} M_{23}+i M_{12} M_{14} M_{24}+i M_{23} M_{24}
 M_{34}) $,

\item 
$M_{12} (M_{12} M_{13} M_{23}+M_{13} M_{14} M_{34}+i M_{14} M_{23} M_{24}) $,

\item 
$M_{12} (M_{12} M_{13} M_{23}+M_{13} M_{14} M_{34}-i M_{14} M_{23} M_{24}) $,

\item 
$ M_{13} (M_{12} M_{13} M_{23}+ M_{13} M_{14} M_{34}+i M_{14} M_{23} M_{24}) $,

\item 
$M_{14} (M_{12} M_{13} M_{23}+ M_{13} M_{14} M_{34}+i M_{14} M_{23} M_{24})$,

\item 
$ M_{14} (M_{12} M_{13} M_{23}+M_{13} M_{14} M_{34}-i M_{14} M_{23} M_{24}) $,

\item 
$ M_{23} (M_{12} M_{13} M_{23}+M_{13} M_{14} M_{34}+i M_{14} M_{23} M_{24}) $,

\item \label{p1}
$ M_{24} (M_{12} M_{13} M_{23}+ M_{13} M_{14} M_{34}+i M_{14} M_{23} M_{24})$,

\item \label{p2}
$ M_{24} (M_{12} M_{13} M_{23}+M_{13} M_{14} M_{34}-i M_{14} M_{23} M_{24}) $,

\item 
$ M_{24} (M_{12} M_{13} M_{23}-M_{13} M_{14} M_{34}+i M_{14} M_{23} M_{24}) $,

\item 
$ M_{34} (M_{12} M_{13} M_{23}+ M_{13} M_{14} M_{34}+i M_{14} M_{23} M_{24}) $,
   
\item 
$ M_{13}^2 M_{23} M_{24}+M_{13} M_{14} M_{23}^2-M_{13} M_{14} M_{34}^2
        +i M_{14} M_{23} M_{24} M_{34} $,

\item 
$ M_{12}^2 M_{13} M_{23}+i M_{12} M_{14} M_{23} M_{24}-M_{13}^2 M_{14} M_{24}        -M_{13} M_{14}^2 M_{23}$,

\item 
$ i\g M_{12} M_{13} M_{23}^2-\g M_{14} M_{23}^2 M_{24}-M_{12} M_{14} M_{24}
   M_{34}-M_{23} M_{24} M_{34}^2 $,

\item 
$ i \g M_{13} M_{14} M_{23} M_{34}-M_{13} M_{14} M_{24}^2-M_{14}^2M_{23}
   M_{24}+M_{23} M_{24} M_{34}^2 $,

\item 
$i \g M_{12} M_{13} M_{14} M_{23}-M_{12}^2 M_{14} M_{24}+
       M_{13} M_{23} M_{24}^2+M_{14} M_{23}^2 M_{24} $,

\item 
$ \g M_{13} M_{14}^2 M_{23}+M_{12} M_{13} M_{23} M_{34}
      +i M_{12} M_{14}^2 M_{24}+M_{13} M_{14} M_{34}^2 $,

\item 
$\g M_{14}^2 M_{23}^2+M_{12}^2M_{14} M_{23}+M_{12} M_{14}^2 M_{34}
      +M_{12} M_{23}^2 M_{34}+M_{14} M_{23} M_{34}^2 $,

\item 
$ -i \g M_{12} M_{13}^2 M_{23}+\g M_{13} M_{14} M_{23} M_{24}
   +M_{12}^2 M_{13} M_{24}+i M_{12} M_{14} M_{24}^2+M_{13} M_{24} M_{34}^2 $,

\item 
$ i \g M_{13}^3 M_{14}+M_{12}^3 M_{13}+i M_{12}^2 M_{14} M_{24}-M_{12} M_{13}
   M_{14}^2+M_{13}^2 M_{24} M_{34} $,

\item 
$\g M_{12} M_{13} M_{14} M_{23}+ M_{12}^2 M_{14} M_{24}
     - M_{12} M_{13} M_{23}^2- M_{13} M_{14} M_{23} M_{34}
     -i M_{13} M_{23} M_{24}^2 $,

\item 
$ i \g M_{13}^2 M_{14} M_{34}+M_{12}^2 M_{13} M_{24}+i M_{12} M_{14} M_{24}^2
       -2 M_{13} M_{14}^2 M_{24}+M_{13} M_{24} M_{34}^2 $,

\item 
$i\g M_{12}^2 M_{13} M_{23}- \g M_{12} M_{14} M_{23} M_{24}
       -i\g M_{13}^2 M_{14} M_{24}+ M_{12} M_{14}^2 M_{24}
       + M_{14} M_{23} M_{24} M_{34} $,

\item 
$ i \g M_{12} M_{13}^2 M_{23}-M_{12}^2M_{13} M_{24}+2 M_{13} M_{23}^2 M_{24}-M_{13}
   M_{24} M_{34}^2+i M_{23} M_{24}^2 M_{34} $,

\item 
$ i \g M_{12}^2 M_{13} M_{23}-M_{12}^3 M_{24}+M_{12} M_{23}^2 M_{24}-
     M_{13} M_{24}^2 M_{34}+i M_{23} M_{24}^3 $,

\item 
$\g M_{14}^2 M_{23} M_{34}-M_{12} M_{14}^2 M_{23}+M_{12} M_{23} M_{34}^2
       -M_{14}^3 M_{34}+i M_{14}^2 M_{24}^2+M_{14} M_{34}^3 $,

\item 
$ i \g M_{13}^3 M_{23}-\g M_{13} M_{14} M_{23} M_{34}-M_{12} M_{13}^2 M_{24}-i M_{12}
   M_{14} M_{24} M_{34}+M_{13} M_{23}^2 M_{34}-M_{13} M_{34}^3 $,

\item 
$\g M_{12} M_{14}^2 M_{23}+i \g M_{13}^2 M_{14}^2+M_{12}^3 M_{14}+
   M_{12}^2 M_{23} M_{34}-M_{12} M_{14}^3-M_{14}^2 M_{23} M_{34}$,

\item 
$ i \g M_{13}^2 M_{23}^2-\g M_{14} M_{23}^2 M_{34}+M_{12} M_{14}
 M_{23}^2-M_{12} M_{14} M_{34}^2+M_{23}^3 M_{34}-M_{23} M_{34}^3 $,

\item 
$ i \g M_{12} M_{14} M_{23}^2+i M_{12}^3 M_{23}+i M_{12}^2 M_{14} M_{34}
      -i M_{12} M_{23}^3-i M_{14} M_{23}^2 M_{34}+M_{23}^2 M_{24}^2 $,

\item 
$ i \g M_{12} M_{13} M_{23} M_{34}-\g M_{14} M_{23} M_{24} M_{34}-
     M_{12} M_{13} M_{24}^2+M_{14}^2 M_{24} M_{34}-i M_{14} M_{24}^3-M_{24} 
     M_{34}^3 $,

\item 
$ i \g M_{12} M_{14} M_{23} M_{34}-i M_{12}^2 M_{14} M_{23}-i M_{12} M_{14}^2
   M_{34}-M_{12} M_{14} M_{24}^2-i M_{12} M_{23}^2 M_{34}-i M_{14} M_{23}
   M_{34}^2+M_{23} M_{24}^2 M_{34} $,

\item 
$i\g M_{12} M_{13}^2 M_{23}- \g M_{12} M_{14} M_{23} M_{34}
    -i \g M_{13}^2 M_{14} M_{34}+ M_{12}^2 M_{14} M_{23}
    + M_{12} M_{14}^2 M_{34}+ M_{12} M_{23}^2 M_{34}
    + M_{14} M_{23} M_{34}^2 $,

\item 
$ \g M_{12}^2 M_{14} M_{23}+i \g M_{12} M_{13}^2 M_{14}+M_{12}^4
        -M_{12}^2 M_{14}^2-M_{12}^2 M_{23}^2-i M_{12} M_{23} M_{24}^2
	+M_{13}^2 M_{24}^2+M_{14}^2 M_{23}^2 $,

\item 
$ -i \g M_{13}^2 M_{23} M_{34}+\g M_{14} M_{23} M_{34}^2+M_{13}^2 M_{24}^2+
     M_{14}^2 M_{23}^2-M_{14}^2 M_{34}^2+i M_{14} M_{24}^2 M_{34}-M_{23}^2
     M_{34}^2+M_{34}^4 $,
\end{enumerate}
\quad\\[-3mm]
where $i^2 = -1$ and $\g \in \kk^\times$.

\bigskip

\begin{ack}
The authors gratefully acknowledge support from the NSF under grants 
DMS-0900239 and DMS-1302050. Moreover, 
the authors are grateful to B.~Shelton for discussions about a potential
approach towards computing the line scheme of the algebra defined 
in~\cite{ST}; that algebra is a member of the family of algebras 
investigated herein.
\end{ack}

%%%%%%%%%%%%%%%%%%%%%%%%%%%%%%%%%%%%%%%%%%%%%%%%%%%%%%%%%%%%%%%%%%%%%
% END  APPENDIX (SEC 5)
% START  BIB
%%%%%%%%%%%%%%%%%%%%%%%%%%%%%%%%%%%%%%%%%%%%%%%%%%%%%%%%%%%%%%%%%%%%%

%\bigskip
\bigskip
\bigskip
%\vfill
%\newpage

\end{document}